\documentclass[12pt,a4paper,reqno]{amsart}
\usepackage{amsmath}
\usepackage{amsfonts}
\usepackage{amssymb}
\numberwithin{equation}{section}

\theoremstyle{plain}

\theoremstyle{definition}

\newsavebox{\proofbox}
\savebox{\proofbox}{\begin{picture}(7,7)%
  \put(0,0){\framebox(7,7){}}\end{picture}}

\def\eps{\varepsilon}

\begin{document}

\title{What is good mathematics?}

\author{Terence Tao}
\address{Department of Mathematics, UCLA, Los Angeles CA 90095-1555, USA.
}
\email{tao@math.ucla.edu}

\begin{abstract} Some personal thoughts and opinions on what ``good quality mathematics'' is, and whether one should try to define this term rigorously.  As a case study, the story of Szemer\'edi's theorem is presented.
\end{abstract}

\maketitle

\section{The many aspects of mathematical quality}

We all agree that mathematicians should strive to produce good mathematics.
But how does one define ``good mathematics'', and should one even dare to try at all?
Let us first consider the former question.  Almost immediately one realises that there are 
many different types of mathematics which could be designated ``good''.  For instance, 
``good mathematics'' could refer (in no particular order) to

\begin{itemize}
\item[(i)] Good mathematical \emph{problem-solving} (e.g. a major breakthrough on an important mathematical problem); 
\item[(ii)] Good mathematical \emph{technique} (e.g. a masterful use of existing methods, or the development of new tools);
\item[(iii)] Good mathematical \emph{theory} (e.g. a conceptual framework or choice of notation which systematically unifies and generalises an existing body of results);
\item[(iv)] Good mathematical \emph{insight} (e.g. a major conceptual simplification, or the realisation of a unifying principle, heuristic, analogy, or theme);
\item[(v)] Good mathematical \emph{discovery} (e.g. the revelation of an unexpected and intriguing new mathematical phenomenon, connection, or counterexample);
\item[(vi)] Good mathematical \emph{application} (e.g. to important problems in physics, engineering, computer science, statistics, etc., or from one field of mathematics to another);
\item[(vii)] Good mathematical \emph{exposition} (e.g. a detailed and informative survey on a timely mathematical topic, or a clear and well-motivated argument);
\item[(viii)] Good mathematical \emph{pedagogy} (e.g. a lecture or writing style which enables others to learn and do mathematics more effectively, or contributions to mathematical education);
\item[(ix)] Good mathematical \emph{vision} (e.g. a long-range and fruitful program or set of conjectures);
\item[(x)] Good mathematical \emph{taste} (e.g. a research goal which is inherently interesting and impacts important topics, themes, or questions);
\item[(xi)] Good mathematical \emph{public relations} (e.g. an effective showcasing of a mathematical achievement to non-mathematicians, or from one field of mathematics to another);
\item[(xii)] Good \emph{meta-mathematics} (e.g. advances in the foundations, philosophy, history, scholarship, or practice of mathematics);
\item[(xiii)] \emph{Rigorous} mathematics (with all details correctly and carefully given in full);
\item[(xiv)] \emph{Beautiful} mathematics (e.g. the amazing identities of Ramanujan; results which are easy (and pretty) to state but not to prove);
\item[(xv)] \emph{Elegant} mathematics (e.g. Paul Erd\H{o}s' concept of ``proofs from the Book''; achieving a difficult result with a minimum of effort);
\item[(xvi)] \emph{Creative} mathematics (e.g. a radically new and original technique, viewpoint, or species of result);
\item[(xvii)] \emph{Useful} mathematics (e.g. a lemma or method which will be used repeatedly in future work on the subject);
\item[(xviii)] \emph{Strong} mathematics (e.g. a sharp result that matches the known counterexamples, or a result which deduces an unexpectedly strong conclusion from a seemingly weak hypothesis);
\item[(xix)] \emph{Deep} mathematics (e.g. a result which is manifestly non-trivial, for instance by capturing a subtle phenomenon beyond the reach of more elementary tools);
\item[(xx)] \emph{Intuitive} mathematics (e.g. an argument which is natural and easily visualisable);
\item[(xxi)] \emph{Definitive} mathematics (e.g. a classification of all objects of a certain type; the final word on a mathematical topic);
\item[(xxii)] etc., etc.\footnote{The above list is not meant to be exhaustive.  In particular, it focuses primarily on the type of mathematics found in mathematical research papers, as opposed to classrooms, textbooks, or papers in disciplines close to mathematics, such as the natural sciences.}
\end{itemize}

As the above list demonstrates, the concept of mathematical quality is a high-dimensional one, and lacks an obvious canonical total ordering\footnote{In particular, it is worth pointing out that mathematical rigour, while highly important, is only one component of what determines a quality piece of mathematics.}.  I believe this is because mathematics is itself complex and high-dimensional, and evolves in unexpected and adaptive ways; each of the above qualities represents a different way in which we as a community improve our understanding and usage of the subject.  There does not appear to be universal agreement as to the relative importance or weight of each of the above qualities.  This is partly due to tactical considerations: a field of mathematics at a given stage of development may be more receptive to one approach to mathematics than another.  It is also partly due to cultural considerations: any given field or school of mathematics tends to attract like-minded mathematicians who prefer similar approaches to a subject.  It also reflects the diversity of mathematical ability; different mathematicians tend to excel in different mathematical styles, and are thus well suited for different types of mathematical challenges.  (See also \cite{gowers-culture} for some related discussion.)

I believe that this diverse and multifaceted nature of ``good mathematics'' is very healthy for mathematics as a whole, as it it allows us to pursue many different approaches to the subject, and 
exploit many different types of mathematical talent, towards our common goal of greater mathematical progress and understanding.   While each one of the above attributes is generally accepted to be a desirable trait to have in mathematics, it can become detrimental to a field to pursue only one or two of them at the expense of all the others.  Consider for instance the following hypothetical (and somewhat exaggerated) scenarios:
\begin{itemize}
\item A field which becomes increasingly ornate and baroque, in which individual results are generalised and refined for their own sake, but the subject as a whole drifts aimlessly without any definite direction or sense of progress; 
\item A field which becomes filled with many astounding conjectures, but with no hope of rigorous progress on any of them; 
\item A field which now consists primarily of using \emph{ad hoc} methods to solve a collection of unrelated problems, which have no unifying theme, connections, or purpose;
\item A field which has become overly dry and theoretical, continually recasting and unifying previous results in increasingly technical formal frameworks, but not generating any exciting new breakthroughs as a consequence; or
\item A field which reveres classical results, and continually presents shorter, simpler, and more elegant proofs of these results, but which does not generate any truly original and new results beyond the classical literature.
\end{itemize}
In each of these cases, the field of mathematics exhibits much activity and progress in the short term, but risks a decline of relevance and a failure to attract younger mathematicians to the subject in the longer term.  Fortunately, it is hard for a field to stagnate in this manner when it is constantly being challenged and revitalised by its connections to other fields of mathematics (or to related sciences), and by exposure to (and respect for) multiple  cultures of ``good mathematics''.  These self-correcting mechanisms help to keep mathematics balanced, unified, productive, and vibrant.  

Let us turn now to the other question posed above, namely whether we should try to pin down a definition of ``good mathematics'' at all.  In doing so, we run the risk of arrogance and hubris; in particular,
we might fail to recognise exotic examples of genuine mathematical progress because they fall outside mainstream definitions\footnote{A related difficulty is that, with the notable exception of mathematical rigour, most of the above qualities are somewhat subjective, and contain some inherent imprecision or uncertainty.  We thank Gil Kalai for emphasising this point.} of ``good mathematics''.  On the other hand, there is a risk also in the opposite position - that all approaches to mathematics are equally suitable and deserving of equal resources\footnote{Examples of scarce resources include money, time, attention, talent, and pages in top journals.} for any given mathematical field of study, or that all contributions to mathematics are equally important; such positions may be admirable for their idealism, but they sap mathematics of its sense of direction and purpose, and can also lead to a sub-optimal allocation of mathematical resources\footnote{Another solution to this problem is to exploit the fact that mathematical resources are also high-dimensional, for instance one can award prizes for exposition, for creativity, etc., or have different journals devoted to different types of achievement.  We thank Gil Kalai for this observation.}.  The true situation lies somewhere in between; for each area of mathematics, the existing body of results, folklore, intuition and experience (or lack thereof) will indicate which types of approaches are likely to be fruitful and thus deserve the majority of resources, and which ones are more speculative and which might warrant inspection by only a handful of independently minded mathematicians, just to cover all bases. For example, in mature and well-developed fields, it may make sense to pursue systematic programs and develop general theories in a rigorous manner, conservatively following tried-and-true methods and established intuition, whereas in newer and less settled fields, a greater emphasis might be placed on making and solving conjectures, experimenting with different approaches, and relying to some extent on non-rigorous heuristics and analogies.  It thus makes sense from a tactical point of view to have at least a partial (but evolving) consensus within each field as to what qualities 
of mathematical progress one should prize the most, so that one can develop and advance the field as effectively as possible at each stage of its development.  For instance, one field may be in great need of solutions to pressing problems; another field may be crying out for a theoretical framework to organise the clutter of existing results, or a grand program or series of conjectures to stimulate new results; other fields would greatly benefit from new, simpler, and more conceptual proofs of key theorems; yet more fields may require good publicity, and lucid introductions to the subject, in order to attract more activity and interest.  Thus the determination of what would constitute good mathematics for a field can and should depend highly on the state of the field itself.  It should also be a determination which is continually updated and debated, both within a field and by external observers to 
that field; as mentioned earlier, it is quite possible for a consensus on how a field should progress to lead 
to imbalances within that field, if they are not detected and corrected in time.

It may seem from the above discussion that the problem of evaluating mathematical quality, while important, is a hopelessly complicated one, especially since many good mathematical achievements may score highly on some of the qualities listed above but not on others; also, many of these qualities are subjective and difficult to measure precisely except with hindsight. However, there is the remarkable phenomenon\footnote{This phenomenon is also somewhat related to the ``unreasonable effectiveness of mathematics'' observed by Wigner \cite{wigner}.} that good mathematics in one of the above senses tends to beget more good mathematics in many of the other senses as well, leading to the tentative conjecture that perhaps there is, after all, a universal notion of good quality mathematics, and all the specific metrics listed above represent different routes to uncover new mathematics, or different stages or aspects of the evolution of a mathematical story.

\section{Case study: Szemer\'edi's theorem}

Turning now from the general to the specific, let us now illustrate the phenomenon mentioned in the preceding paragraph by considering the history and context of \emph{Szemer\'edi's theorem} \cite{szemeredi} - the beautiful and celebrated result that any subset of integers of positive (upper) density must necessarily contain arbitrarily long arithmetic progressions.   
I will avoid all technical details here; the interested reader is referred to \cite{dichotomy} and the references therein for further discussion.

There are several natural places to start this story.  I will begin with \emph{Ramsey's theorem} \cite{ramsey}:
that any finitely coloured, sufficiently large complete graph will contain large monochromatic complete subgraphs. (For instance, given any six people, either three will know each other, or three will be strangers to each other, assuming of course that ``knowing one another'' is a well-defined and symmetric relation.)  This result, while simple to prove (relying on nothing more than an iterated pigeonhole principle), represented the discovery of a new phenomenon and created a new species of mathematical result: the \emph{Ramsey-type theorem}, each one of which being a different formalisation of the newly gained insight in mathematics that \emph{complete disorder is impossible}. 

One of the first Ramsey-type theorems (which actually predates Ramsey's theorem by a few years) was \emph{van der Waerden's theorem} \cite{vdw}: given any finite colouring of the integers, one of the colour classes must contain arbitrarily long arithmetic progressions.  Van der Waerden's highly recursive proof was very elegant, but had the drawback that it offered fantastically poor \emph{quantitative} bounds for the appearance of the first arithmetic progression of a given length; indeed, the bound involved an Ackermann function of this length and the number of colours.  Erd\H{o}s and Tur\'an \cite{erdos} had the good mathematical taste to pursue this quantitative question\footnote{Erd\H{o}s also pursued the question of quantitative bounds for the original theorem of Ramsey, leading among other things to the founding of the immensely important \emph{probabilistic method} in combinatorics, but this is a whole story in itself which we have no space to discuss here.} further, being motivated also by the desire to make progress on the (then conjectural) problem of whether the primes contained arbitrarily long progressions.  They then advanced a number of strong conjectures, one of which became Szemer\'edi's theorem; another was the beautiful but (still open) stronger statement that any set of positive integers whose reciprocals were not absolutely summable contained arbitrarily long arithmetic progressions.   

The first progress on these conjectures was a sequence of counterexamples, culminating in the elegant construction of
Behrend \cite{behrend} of a moderately sparse set (whose density in $\{1,\ldots,N\}$ was asymptotically greater than $N^{-\eps}$ for any fixed $\eps$) without arithmetic progressions of length three.  This construction ruled out the most ambitious of the Erd\H{o}s-Tur\'an conjectures (in which polynomially sparse sets were conjectured to have many progressions), and as a consequence also ruled out a significant class of elementary approaches to these problems (e.g. those based on inequalities such as the Cauchy-Schwarz or H\"older inequalities).  While these examples did not fully settle the problem, they did indicate that the Erd\H{o}s-Tur\'an conjectures, if true, would necessarily have a non-trivial (and thus presumably interesting) proof.

The next major advance was by Roth \cite{roth}, who applied the \emph{Hardy-Littlewood circle method}\footnote{Again, the history of the circle method is another great story which we cannot detail here.  Suffice to say though that this method, in modern language, is part of the now standard insight that Fourier analysis is an important tool for tackling problems in additive combinatorics.} together with a new method (the \emph{density increment argument}) in a beautifully elegant manner to establish \emph{Roth's theorem}: every set of integers of positive density contained infinitely many progressions of length three.  It was then natural to try to extend Roth's methods to progressions of longer length. Roth and many others tried to do so for many years, but without full success; the reason for the obstruction here was not fully appreciated until the work of Gowers much later.  It took the formidable genius of Endr\'e Szemer\'edi \cite{szemeredi-4}, \cite{szemeredi}, who returned to purely combinatorial methods (in particular, pushing the density increment argument to remarkable new levels of technical sophistication) to extend Roth's theorem first to progressions of length four\footnote{Shortly afterward, Roth \cite{roth-4} was able to combine some of Szemer\'edi's ideas with his own Fourier analytic method to create a hybrid proof of Szemer\'edi's theorem for progressions of length four.}, and then to progressions of arbitrary length, thus establishing his famous theorem.  Szemer\'edi's proof was a technical \emph{tour de force}, and introduced many new ideas and techniques, the most important of which was a new way to look at extremely large graphs, namely to approximate them by bounded complexity models.  This result, the celebrated and very useful \emph{Szemer\'edi regularity lemma}, is notable on many levels.  As mentioned above, it gave a radically new insight regarding the structure of large graphs (which in modern language is now regarded as a \emph{structure theorem} as well as a  \emph{compactness theorem} for such graphs); it gave a new proof method (the \emph{energy increment method}) which will become crucial later in this story; and it also generated an incredibly large number of unexpected applications, from graph theory to property testing to additive combinatorics; the full story of this regularity lemma is unfortunately too lengthy to be described here.

While Szemer\'edi's accomplishment is undoubtedly a highlight of this particular story, it was by no means the last word on the matter.  Szemer\'edi's proof of his theorem, while elementary, was remarkably intricate, and not easily comprehended.  It also did not fully resolve the original questions motivating Erd\H{o}s and Tur\'an, as the proof itself used van der Waerden's theorem at two key junctures and so did not give any improved quantitative bound on that theorem.  Furstenberg then had the mathematical taste to seek out a radically different (and highly non-elementary\footnote{For instance, some versions of Furstenberg's argument rely heavily on the axiom of choice, though it is possible to also recast the argument in a choice-free manner.}) proof, based on an insightful analogy between combinatorial number theory and ergodic theory which he soon formalised as the very useful \emph{Furstenberg correspondence principle}.  From this principle\footnote{There is also a similar correspondence principle which identifies van der Waerden's theorem with a multiple recurrence theorem for \emph{topological} dynamical systems.  This leads to the fascinating story of \emph{topological dynamics}, which we unfortunately have no space to describe here.} one readily concludes that Szemer\'edi's theorem is equivalent to a \emph{multiple recurrence theorem} for measure-preserving systems.  It then became natural to prove this theorem (now known as the \emph{Furstenberg recurrence theorem}) directly by methods from ergodic theory, in particular by exploiting the various classifications and structural decompositions (e.g. the ergodic decomposition) available for such systems.  Indeed, Furstenberg soon established the \emph{Furstenberg structure theorem}, which described any measure preserving system as a weakly mixing extension of a tower of compact extensions of a trivial system, and based on this theorem and several additional arguments (including a variant of the van der Waerden argument) was able to establish the multiple recurrence theorem, and thus give a new proof of Szemer\'edi's theorem.  It is also worth mentioning that Furstenberg also produced an excellent book \cite{furst-book} on this and related topics, which systematically formalised the basic theory while also contributing greatly to the growth and further development of this area.

Furstenberg and his coauthors then realised that this new method was potentially very powerful, and could be used to establish many more types of recurrence theorems, which (via the correspondence principle) then would yield a number of highly non-trivial combinatorial theorems. Pursuing this vision, Furstenberg, Katznelson, and others obtained many variants and generalisations of Szemer\'edi's theorem, obtaining for instance variants in higher dimensions and even establishing a density version of the Hales-Jewett theorem \cite{hales-jewett} (a very powerful and abstract generalisation of the van der Waerden theorem).  Many of the results obtained by these infinitary ergodic theory techniques are not known, even today, to have any ``elementary'' proof, thus testifying to the power of this method.  Furthermore, as a valuable byproduct of these efforts, a much deeper understanding of the structural classification of measure-preserving systems was obtained.  In particular, it was realised that for many classes of recurrence problem, the asymptotic recurrence properties of an arbitrary system are almost completely controlled by a special \emph{factor} of that system, known as the (minimal) \emph{characteristic factor} of that system\footnote{An early example of this is von Neumann's mean ergodic theorem, in which the factor of shift-invariant functions controls the limiting behaviour of simple averages of shifts.}.  
Determining the precise nature of this characteristic factor for various types of recurrence then became a major focus of study, as it was realised that this would lead to more precise information on the limiting behaviour (in particular, it would show that certain asymptotic expressions related to multiple recurrence actually converged to a limit, which was a question left open from Furstenberg's original arguments).  Counterexamples of Furstenberg and Weiss, as well as results of Conze and Lesigne, eventually led to the conclusion that these characteristic factors should be describable by a very special (and algebraic) type of measure-preserving system, namely a \emph{nilsystem} associated with nilpotent groups; these conclusions culminated in precise and rigorous descriptions of these factors in a technically impressive paper of Host and Kra \cite{host-kra2} (and subsequently also by Ziegler \cite{ziegler}), which among other things settled the question mentioned earlier concerning convergence of the asymptotic multiple recurrence averages.  The central role of these characteristic factors illustrated quite starkly the presence of a dichotomy between structure (as represented here by nilsystems), and randomness (which is captured by a certain technical type of ``mixing''
property), and to the insight that it is this dichotomy which in fact underlies and powers Szemer\'edi's theorem and
its relatives.  Another feature of the Host-Kra analysis worth mentioning is the prominent appearance of averages associated to ``cubes'' or ``parallelopipeds'', which turn out to be more tractable to analyse for a number of reasons than the multiple recurrence averages associated to arithmetic progressions.  

In parallel to these ergodic theory developments, other mathematicians were seeking to understand, reprove, and improve upon Szemer\'edi's theorem in other ways.  An important conceptual breakthrough was made by Ruzsa and Szemer\'edi \cite{rsz}, who used the Szemer\'edi regularity lemma mentioned earlier to establish a number of results in graph theory, including what is now known as the \emph{triangle removal lemma}, which roughly asserts that a graph which contains a small number triangles can have those triangles removed by deleting a surprisingly small number of edges.  They then observed that the Behrend example mentioned earlier gave some limits as to the quantitative bounds in this lemma, in particular ruling out many classes of elementary approaches to this lemma (as such approaches typically give polynomial type bounds); indeed to this day all known proofs of the removal lemma proceed via some variant of the regularity lemma.  Applying this connection in the contrapositive, it was observed that in fact the triangle removal lemma implied Roth's theorem on progressions of length three.  This discovery opened up for the first time 
the possibility that Szemer\'edi type theorems could be proven by purely \emph{graph-theoretical} techniques, discarding almost entirely the additive structure of the problem.  (Note that the ergodic theory approach still retained this structure, in the guise of the action of the shift operator on the system; also, Szemer\'edi's original proof is only partly graph-theoretical, as it exploits the additive structure of progressions in many different places.)  It took some time though to realise that the graph theoretic method, like the Fourier-analytic method before it, was largely restricted to detecting ``low complexity'' patterns such as triangles or progressions of length three, and to detect progressions of longer length would require the substantially more difficult theory of \emph{hypergraphs}.  In particular this motivated the program (spearheaded by Frankl and R\"odl) for obtaining satisfactory analogue of the regularity lemma for hypergraphs, which would be strong enough to yield consequences such as Szemer\'edi's theorem (as well as a number of variants and generalisations).  This turned out to be a surprisingly delicate task, in particular carefully arranging the hierarchy\footnote{This hierarchy seems related to the towers of extensions encountered by Furstenberg in his analogous quest to ``regularise'' a measure-preserving system, though the precise connection is still poorly understood at present.} of parameters involved in such a regularisation so that they dominated each other in the correct order.  Indeed the final versions of the regularity lemma, and the companion ``counting lemmas'' from which one could deduce Szemer\'edi's theorem, have only appeared rather recently (\cite{nrs}, \cite{rs}, \cite{rodl}, \cite{rodl2}, \cite{gowers-reg}, $\ldots$).  It is also worth mentioning a very instructive counterexample \cite{gowers-sz} of Gowers, which shows that the quantitative bounds in the original regularity lemma must be at least tower-exponential in nature, thus indicating again the non-trivial nature (and power) of this lemma.

The Fourier analytic approach to Szemer\'edi's theorem, which had not progressed significantly since the work of Roth, was finally revisited by Gowers \cite{gowers-4}, \cite{gowers}.  As with other approaches, the Fourier-analytic approach proceeded by establishing a dichotomy on sets of integers, that they were either \emph{structured} or \emph{pseudorandom} in some sense.  The relevant notion of structure here was worked out by Roth - structured sets should enjoy a density increment on medium-length arithmetic progressions - but the correct notion of pseudorandomness or ``uniformity'' was less clear.  Gowers produced an example (closely related, in fact, to examples of Furstenberg and Weiss mentioned earlier) showing that Fourier-based notions of pseudorandomness were inadequate for controlling progressions of length four and higher, and then proceeded to introduce a different notion of uniformity (very closely related to the cube averages of Host and Kra, and also to certain notions of hypergraph regularity) which sufficed.  The remaining task was to establish a quantitative and rigorous form of the dichotomy.  This turned out to be
surprisingly difficult (mainly due to the limited utility of the Fourier transform in this setting), and in many ways
analogous to the efforts of Host-Kra and Ziegler to endow characteristic factors with the algebraic structure of
nilsystems.  However, by combining Fourier analytic tools with major results from additive combinatorics such as Freiman's theorem and the Balog-Szemer\'edi theorem (the history of these being also an interesting story in its own right, see e.g. \cite{tao-vu}), together with several new combinatorial and probabilistic methods, Gowers was able to achieve this in a remarkable \emph{tour de force}, and in particular obtained remarkably strong quantitative bounds on Szemer\'edi's theorem and van der Waerden's theorem\footnote{It is worth noting also the brilliantly creative proof of van der Waerden's theorem by Shelah \cite{shelah}, which held the previous record for the best constants for this theorem; the ideas of Shelah's proof have not yet been successfully integrated into the rest of the subject, but I expect that this will happen in the future.}.

To summarise so far, four parallel proofs of Szemer\'edi's theorem have been achieved; one by direct combinatorics, one by ergodic theory, one by hypergraph theory, and one by Fourier analysis and additive combinatorics.  Even with so many proofs, there was still a sense that our understanding of this result was still incomplete; for instance, none of the approaches were powerful enough to detect progressions in the primes, mainly because of the sparsity of the prime sequence.  (The Fourier method, or more precisely the Hardy-Littlewood-Vinogradov circle method, can be used however to establish infinitely many progressions of length three in the primes \cite{van-der-corput}, and with substantially more effort can also partially address progressions of length four \cite{heath-brown1}.)  However, by using ideas from restriction theory in harmonic analysis (which is another fascinating story that we will not discuss here), Green \cite{green} was able to treat the primes ``as if'' they were dense, and in particular obtain an analogue of Roth's theorem for dense subsets of primes.  This opened up the intruiging possibility of a \emph{relative Szemer\'edi theorem}, allowing one to detect arithmetic progressions in dense subsets of other sets than the integers, for instance dense subsets of primes.  Indeed, a prototypical relative Roth theorem for dense subsets of quite sparse random sets had already
appeared in the graph theory literature \cite{klr}.

In joint work\footnote{Incidentally, I was initially attracted to these problems by their intersection with another great mathematical story, that of the \emph{Kakeya conjecture}, which we again do not have space to discuss here.  It is however related in a somewhat surprising fashion with the story of restriction theory mentioned earlier.} with Ben Green, we began the task of trying to relativise Gowers' Fourier analytic and combinatorial arguments to such contexts as dense subsets of sparse random or ``pseudorandom'' sets.  After much effort (inspired in part by the hypergraph theory, which was well adapted to count patterns in sparse sets, and also in part by an ``arithmetic regularity lemma'' of Green \cite{gr2} that adapted the regularity lemma ideas from graph theory to additive contexts)
we were eventually able (in an unpublished work) to detect progressions of length four in such sets.  At this point
we realised the analogies between the regularity lemma approach we were using, and the characteristic factor
constructions in Host-Kra, and by substituting\footnote{This was a little tricky for a number of reasons, most notably that the ergodic theory constructions were infinitary in nature, whereas to deal with the primes it was necessary to work in a finitary context.  Fortunately I had already attempted to finitise the ergodic approach to Szemer\'edi's theorem \cite{tao:ergodic}; while that attempt was incomplete at the time, it turned out to have enough substance to be helpful for the application to the primes.} in those constructions (in particular relying heavily on cube averages) we were able to establish a satisfactory relative Szemer\'edi theorem, which relied on a certain \emph{transference principle} which asserted, roughly speaking, that dense subsets of sparse pseudorandom sets behaved ``as if'' they were dense in the original space.  In order to apply this theorem to the primes, we needed to enclose the primes in a suitably pseudorandom set (or more precisely a measure); but very fortuitiously for us, the recent breakthroughs\footnote{At the time of our paper, the construction we used was from a paper of Goldston and Y{\i}ld{\i}r{\i}m that was retracted for an unrelated error, which they eventually repaired with some clever new ideas from Pintz \cite{pintz}.  This supports a point made earlier, that it is not absolutely necessary for a piece of mathematics to be absolutely correct in every detail in order to be of value to future (rigorous) work.} of Goldston and Y{\i}ld{\i}r{\i}m \cite{goldston-yildirim} on prime gaps\footnote{Once again, the story of prime gaps is an interesting one which we will be unable to pursue here.} had constructed almost exactly what we needed for this purpose, allowing us to establish at last the old conjecture that the primes contained arbitrarily long arithmetic progressions.

The story still does not end here, but instead continues to develop in several directions.  On one hand there are now a number of further applications of the transference principle, for instance to obtain constellations in the Gaussian primes or polynomial progressions in the rational primes.  Another promising avenue of research is the convergence of the Fourier, hypergraph, and ergodic theory methods to each other, for instance in developing infinitary versions of the graph and hypergraph theory (which have applications to other areas of mathematics as well, such as property testing) or finitary versions of the ergodic theory.  A third direction is to make the nilsystems that control recurrence in the ergodic theory setting, also control various finitary averages relating to arithmetic progressions; in particular, Green and I are actively working on computing correlations between primes and sequences generated by nilsystems (using methods dating back to Vinogradov) in order to establish more precise asymptotics on various patterns that can be found in the primes.  Last, but not least, there is the original Erd\H{o}s-Tur\'an conjecture, which still remains open after all this progress, though there is now some very promising advances of Bourgain \cite{bourgain-triples}, \cite{bourgain-preprint} which should lead to further developments.

\section{Conclusion}

As we can see from the above case study, the very best examples of good mathematics do not merely fulfil one or more of the criteria of mathematical quality listed at the beginning of this article, but are more importantly part of a greater mathematical \emph{story}, which then unfurls to generate many further pieces of good mathematics of many different types.  Indeed, one can view the history of entire fields of mathematics as being primarily generated by a handful of these great stories, their evolution through time, and their interaction with each other.  I would thus conclude that good mathematics is not merely measured by one or more of the ``local'' qualities listed previously (though these are certainly important, and worth pursuing and debating), but also depends on the more ``global'' question of how it fits in with other pieces of good mathematics, either by building upon earlier achievements or encouraging the development of future breakthroughs.  Of course, without the benefit of hindsight it is difficult to predict with certainty what types of mathematics will have such a property.  
There does however seem to be some undefinable sense that a certain piece of mathematics is ``on to something'', that it is a piece of a larger puzzle waiting to be explored further.  And it seems to me that it is the pursuit of such 
intangible promises of future potential are least as important an aspect of mathematical progress than the more
concrete and obvious aspects of mathematical quality listed previously.  Thus I believe that good mathematics is more than simply the process of solving problems, building theories, and making arguments shorter, stronger, clearer, more elegant, or more rigorous, though these are of course all admirable goals; while achieving all of these tasks (and debating which ones should have higher priority within any given field), we should also be aware of any possible larger context that one's results could be placed in, as this may well lead to the greatest long-term benefit for the result, for the field, and for mathematics as a whole.

\section{Acknowledgements}

I thank Laura Kim for reading and commenting on an early draft of this manuscript, and Gil Kalai for many thoughtful comments and suggestions.

\end{document}